\documentclass[12pt]{amsart}
\usepackage[cp1251]{inputenc}

\usepackage{amssymb,amsmath,amscd,color}

\theoremstyle{plain}
\newtheorem{theorem}{Theorem}
\newtheorem{proposition}{Proposition}

\newcommand{\cA}{{\mathcal{A}}}
\newcommand{\cD}{{\mathcal{D}}}
\newcommand{\cP}{{\mathcal{P}}}
\newcommand{\cU}{{\mathcal{U}}}

        \newcommand{\field}[1]{{\mathbb{#1}}}
        \newcommand{\zz}{\field{Z}}
        \newcommand{\rr}{\field{R}}
        \newcommand{\cc}{\field{C}}

\allowdisplaybreaks

\begin{document}

\title[]{Floquet--Bloch functions on non-simply connected manifolds, the Aharonov--Bohm fluxes, and conformal invariants of immersed surfaces}
\author[]{I.A. Taimanov}
\address{Sobolev Institute of Mathematics, Novosibirsk, 630090 Russia} 
\email{taimanov@math.nsc.ru}

\begin{abstract}
Spectral (Bloch) varieties of multidimensional diffe\-ren\-tial operators on non-simply connected manifolds are defined. In their terms it is given a description of the analytical dependence of the spectra of magnetic Laplacians on non-simply connected manifolds on the values of the Aharonov--Bohm fluxes and a
con\-struc\-tion of analogues of spectral curves for two-dimensional Dirac operators on Riemann surfaces 
and, thereby, new conformal inva\-riants of immersions of surfaces into $3$- and $4$-dimensional Euclidean spaces.
\end{abstract}

\date{}

\maketitle

\hfill{\sl To S.P. Novikov on his 85-th birthday}

\vskip5mm

\section{Introduction}

The modern theory of finite-zone (finite gap) operators goes back to Novikov's work on the integration of the periodic Korteweg--de Vries equation \cite{N}. The spectral curves of finite-zone one-dimensional opera\-tors introduced in it Schr\"odinger later formed the basis of algebraic-geometric methods for integrating soliton equations. The generalization of spectral curves to the case of two-dimensional Schr\"odinger operators, defined in \cite{DKN}, has expanded the applications of this method.

In this paper we will discuss applications of spectral (Bloch) varieties of multidimensional differential operators on non-simply connected ma\-ni\-folds.

The main objectives of this work are

\begin{itemize}
\item
description of the analytical dependence of the spectra of magne\-tic Laplacians on non-simply connected 
ma\-ni\-folds on the values of the Aharonov--Bohm fluxes (\S 3, Theorem 2);

\item
construction of analogues of spectral curves for two-dimensional Dirac ope\-ra\-tors on Riemann surfaces and, thereby, new confor\-mal invariants of immersions of surfaces
into $3$- and $4$-dimensio\-nal Euclidean spaces (\SS 4 and 5, Theorems 4 and 5).
\end{itemize}

Let us point out as an example another nontrivial situation, when the spectral curve (one-dimensional Bloch variety) is naturally defined. In \cite{BT} there were considered
Bloch eigenfunctions of the Cauchy--Riemann operator on a punctured torus
$$
T^2 = \cc / \{\zz e_1 + \zz e_2\},
$$
i.e. a two-dimensional torus with punctured points $p_1,\dots,p_N$. They are defined as solutions to the equation
$$
\bar{\partial}\psi = 0,
$$
satisfying the Floquet--Bloch conditions
$$
\psi(z + e_j) = e^{2\pi i \langle k,e_j\rangle} \psi(z),
$$
and having asymptotics
$$
\psi(z) = \frac{a_j}{z - p_j} + O(|z-p_j|.
$$
It was shown in \cite{BT} that such functions
exist exactly for tuples
$$
(k_1,k_2,a_1,\dots,a_N) \in \Gamma,
$$
belonging to some spectral curve $\Gamma$.

After our talk at a conference in Lumumba Russian People's Friend\-ship University in October 2023, Yu. Kordyukov drew our attention to the works of \cite{Kha,KOS,Sunada} devoted to the Floquet--Bloch theory for periodic elliptic operators on non-simply connected compact manifolds. They contain part of the results from \S 2, but we save their conclusion, since it was not presented separately and will be needed later.
These works study a generalization of the Bloch property and
its application to questions of the existence of eigenvalues of such operators. We will be more interested in Bloch varieties.

\section{Floquet--Bloch eigenfunctions of Laplace--Beltrami and Schr\"odinger operators}

Let $M$ be a closed manifold with a nontrivial fundamental group $\pi_1(M)$.
Consider the composition epimorphism
$$
\rho: \pi_1(M) \to H_1(M) \to \zz^k = H_1(M)/\mathrm{Torsion}
$$
and denote by
$$
p: M_0 \to M
$$
the covering corresponding to the subgroup $\mathrm{Ker}\,\rho$:
$\pi_1(M_0) = \mathrm{Ker}\, \rho.$
The group 
$\Lambda = H_1(M)/\mathrm{Torsion} = \zz^N$, where $N = \dim H_1(M;\rr)$
acts freely on $M_0$, so that $M = M_0/\Lambda$.

Let us choose in $\Lambda$ the generators $Z_1,\dots,Z_N$ and in $H^1(M;\rr)$ a basis 
realized by closed $1$-forms $\omega_1,\dots,\omega_N$, such that
$$
\langle \omega_j, Z_k \rangle = \int_{Z_j} \omega_k = \delta_{jk}, \ \ \ j,k=1,\dots,N.
$$
The preimages under covering $p$ of $1$-forms $\omega_1,\dots,\omega_N$
are the differentials of functions
$$
p^\ast \omega_k = dh_k, \ \ h_k: M_0 \to \rr, \ \ k=1,\dots,N.
$$
Let us denote by $T_1,\dots,T_N$ generators of the action of $\Lambda$ on $M_0$ such that
$$
h_j (T_k x) = h_j(x) + \delta_{jk}.
$$

Let us consider functions on $M_0$ satisfying the conditions
\begin{equation}
\label{fb}
\begin{split}
\psi: M_0 \to \cc, \\
\psi(T_k x) = e^{i \kappa_m} \psi(x), \ \ \ m=1,\dots,N,
\end{split}
\end{equation}
in this case the quantities $\kappa_1,\dots,\kappa_N \in \cc$ are called
{\it quasi-momenta} and are determined up to $2\pi m, m \in \zz$.
\footnote{If the periods of the functions $h_m$ are not normalized by the unit, then the quasimomenta are determined
taking into account these periods. for example, in the theory of the one-dimensional periodic Schr\"odinger operator
after a shift by period $T$, the Bloch function is multiplied by $e^{i \kappa T}$, where $\kappa$ is the quasi-momentum \cite{N}.}

The multipliers which correspond to such a function define a homo\-mor\-phism
$$
\mu: \Lambda \to \cc^\ast = \cc \setminus \{0\}, \ \ \mu(T_m) = e^{i \kappa_m}, \ m=1,\dots,N.
$$

We have

\begin{proposition}
\label{p1}
If a function $\psi$ satisfies \eqref{fb}, then it can be represented as
\begin{equation}
\label{b}
\psi(x) = e^{i(\kappa_1 h_1(x) + \dots + \kappa_N h_N(x))} \varphi(x),
\end{equation}
where $\varphi: M_0 \to \cc$ is a $\Lambda$-periodic function.
\end{proposition}

The Floquet-Bloch functions were defined as
solutions of differential equations with periodic coefficients on $\rr^N$,
satisfying \eqref{fb}, Floquet for $N=1$ and Bloch for arbitrary dimension $N$.
It was assumed that they have the form \eqref{b}, where
$h_1(x) = x^1,\dots, h_N(x) = x^N$ --- coordinate functions.

Let us consider the Laplace operator on the Riemannian manifold $M$. Let $g_{jk}dx^j dx^k$ denote the Riemannian metric. The Laplace--Beltrami operator (or, more simply, the Laplacian), corresponding to this metric, acts on smooth functions on the manifold and is given by the expression
\begin{equation}
\label{laplace}
\Delta = \frac{1}{\sqrt{g}}\, \frac{\partial}{\partial x^j}\, \sqrt{g} \, g^{jk}\, \frac{\partial}{\partial x^k},
\end{equation}
where $g = \det (g_{jk})$, $g^{jk}g_{km} = \delta^j_m$ and everywhere we
mean summation over repeated upper and lower indices. Its action on the function $f$ takes the form
$$
\Delta f = \frac{1}{\sqrt{g}}\, \frac{\partial}{\partial x^j}\left( \sqrt{g} \, g^{jk}\, \frac {\partial f}{\partial x^k}\right).
$$
The Schr\"odinger operator on a manifold is defined as the sum of the operator $-\Delta$ and multiplication by the potential $U(x)$:
$$
L = -\Delta + U.
$$

We call by {\it Bloch eigenfunction} of the Schr\"odinger operator with eigenvalue $E$ a
solution to the equation
\begin{equation}
\label{eigenf1}
L\psi = E\psi,
\end{equation}
satisfying \eqref{fb}. As a function it is defined on
the maximal Abelian covering $M_0$ of $M = M_0/\zz^N$.

The following fact is proved by straightforward  calculations.

\begin{proposition}
\label{p2}
The function $\psi: M_0 \to \cc$ of the form \eqref{b} satisfies the equation
\begin{equation}
\label{sch}
L\psi = (-\Delta + U)\psi = E\psi
\end{equation}
if and only if the function $\varphi: M \to \cc$
satisfies the equation
\begin{equation}
\label{eqphi}
\begin{split}
-\left[\Delta \varphi - g^{jk} \frac{\partial}{\partial x^j} \left(\sum_{m=1}^N \kappa_m h_m\right) \frac{\partial }{\partial x^k}
\left(\sum_{q=1}^N \kappa_q h_q \right)
\varphi + \right.
\\
\left. + i \left(\sum_{m=1}^N \kappa_m \Delta h_m \right) \varphi
  + 2i g^{jk} \left(\sum_{m=1}^N \kappa_m \frac{\partial h_m}{\partial x^j}\right) \frac{\partial \varphi}{\partial x ^k}
\right]
+ U\varphi = E\varphi.
\end{split}
\end{equation}
\end{proposition}

It should be noted that the terms in \eqref{eqphi} that include the functions $h_m$ are correctly defined on $M$, although the functions $h_m$ are defined on the cover $M_0$. By definition, the function $h_m$ is chosen such that
$$
dh_m = p^\ast \omega_m,
$$
where $\omega_m = A_{m,k}dx^k$ is a smooth $1$-form on the compact manifold $M$.
That's why
$$
\frac{\partial h_m}{\partial x^k} = A_{m,k}, \ \ \ m=1,\dots,N, k=1,\dots, \dim M,
$$
and this quantity is correctly defined on $M$. Let's write the equation \eqref{eqphi}
in terms of $1$-forms $A$:
\begin{equation}
\label{eqphi1}
\begin{split}
-\left[\Delta - g^{jk} \left(\sum_{m=1}^N \kappa_m A_{m,j} \right)
\left(\sum_{q=1}^N \kappa_q A_{q,k} \right)  + \right.
\\
+ i \sum_{m=1}^N \kappa_m \left( g^{jk}\frac{\partial A_{m,k}}{\partial x^j} + \frac{1}{\sqrt{g}}
\frac{\partial \sqrt{g}g^{jk}}{\partial x^j} A_{m,k}\right) +
\\
\left. + 2i g^{jk} \left(\sum_{m=1}^N \kappa_m A_{m,j}\right) \frac{\partial }{\partial x^k}
\right]  \varphi
+ U \varphi = E \varphi.
\end{split}
\end{equation}

Recall that $L_p(M), p\in \zz, p \geq 0,$ denotes the space of functions $M$ such that
$\int_M |f|^p d\mu < \infty$,
where $d\mu$ is the measure generated by the metric.
More precisely, this is the space of equivalence classes of functions with respect to the relation
$f_1 \sim f_2$ if $f_1$ and $f_2$ coincide outside a set of measure zero.
In $L_p(M)$ the norm $\|f\|_p = \left(\int_M |f|^p d\mu\right)^{\frac{1}{p}}$ is defined.
For $C^\infty$-functions on $M$ we define norms
$$
\|\nabla^k f\|^2 = \nabla^{j_1}\dots\nabla^{j_k} f \nabla_{j_1}\dots\nabla_{j_k} f,
$$
where $\nabla$ is the covariant derivative.
We denote by $H^p_k(M)$ the closure of the space of $C^\infty$-functions
with $\|\nabla^l f\| \in L_p(M)$ for $0 \leq l \leq k$ relative to the norm
$\|g\|_{H^p_k} = \sum_{l=0}^k \|\nabla^l f\|_p$.
If the Riemannian manifold $M$ is compact and without boundary, then the spaces defined in this way do not depend on the choice of metric and Kondrashov’s theorem holds:
the embedding $H^q_k \subset L_p$ is compact for $1 \geq \frac{1}{p} > \frac{1}{q} - \frac{k}{\dim M}$.
Detailed definitions and proofs of these facts are presented in \cite{Aubin}.

\begin{theorem}
In $\cc^{N+1}$ there is an analytic subset $\Gamma$, a spectral variety distinguished as the set
zeros
$$
\Gamma = \{\cP(E,\kappa_1,\dots,\kappa_N)=0\}
$$
of a nontrivial entire function such that a point $(E,\kappa_1,\dots,\kappa_N)$ belongs to $\Gamma$ if and only if there exists a Bloch eigenfunction of the form \eqref{b}
with quasimomenta $\kappa_1,\dots,\kappa_N$, satisfying \eqref{sch}.

The spectral variety does not depend on the choice of functions $h_1$, $\dots$, $h_N$ such that
$p^\ast \omega_k = dh_k, k=1,\dots,N$, and therefore is determined only by the operator $L$.
\end{theorem}

{\sc Proof.} Let us choose a constant $\varepsilon$ such that the operator 
$\Delta + \varepsilon: H^2_2(M) \to L_2(M)$
is invertible and represent $\varphi$ in the form
$$
\varphi = (\Delta+\varepsilon)^{-1}\xi.
$$
We rewrite the equation \eqref{eqphi1} in the form
\begin{equation}
\label{eqphi2}
\begin{split}
\left\{1 + \left[(E-U-\varepsilon) + \sum_m \kappa_m B_m + \sum_{m,q} \kappa_m \kappa_q B_{mq} +
\right.\right.
\\
\left.\left.
+ \sum_m \kappa_m C_m^k \frac{\partial}{\partial x^k}\right](\Delta+\varepsilon)^{-1}\right\}\xi = 0.
\end{split}
\end{equation}
Composition of the multiplication by a bounded function
$G_\kappa = [(E-U-\varepsilon) + \sum_m \kappa_m B_m + \sum_{m,q} \kappa_m \kappa_q B_{mq}]$
and $(\Delta+\varepsilon)^{-1}$ splits into a composition
$$
L_2 \stackrel{(\Delta+\varepsilon)^{-1}}{\longrightarrow} H^2_2 \stackrel{\mbox{embedding}}{\longrightarrow}
\ \ L_2 \stackrel{\times G_\kappa}{\longrightarrow} L_2.
$$
The composition of the first two operators is an operator that increases smoothness, 
the first being bounded, and the embedding is compact. The multiplication operator by $G_\kappa$ is bounded.
The composition of bounded and compact operators is compact, hence the operator
$$
[(E-U-\varepsilon) + \sum_m \kappa_m B_m + \sum_{m,q} \kappa_m \kappa_q B_{mq}](\Delta+\varepsilon)^{-1}
$$
is compact for any $\kappa$ and $E$.
There is a similar decomposition
$$
L_2 \stackrel{(\Delta+\varepsilon)^{-1}}{\longrightarrow} H^2_2 \stackrel{\sum \kappa_m C^k_m \frac{\partial}{\partial x^k}}{\longrightarrow}
\ \ H^1_2 \stackrel{\mbox{embedding}}{\longrightarrow} L_2,
$$
in which the embedding is compact and other operators are bounded.

Therefore, the equation \eqref{eqphi2} has the form
\begin{equation}
\label{kel}
(1+ A(E,\kappa_1,\dots,\kappa_N))\xi = 0,
\end{equation}
where $A(E,\kappa_1,\dots,\kappa_N)$ is a polynomial in $\kappa_1,\dots,\kappa_N$ and $E$ is a pencil of compact operators from $L_2(M)$ to $L_2(M )$.
By the Keldysh theorem, in this case the equation \eqref{kel} is solvable
then and only if
$(E,\kappa_1,\dot,\kappa_N)$ belongs to some analytic subset $\Gamma$
in $\cc^{N+1}$,
which is the zero set of an entire function $\cP$.

Since for $\kappa_1=\dots=\kappa_N = 0$ this subset intersects $\cc$ along the spectrum of the operator $-\Delta+E$, which is discrete, then $\Gamma$ is a proper subset and the entire function 
$\cP$ is non-trivial.

Suppose we took other functions $\tilde{h}_1,\dots,\tilde{h}_N$,
satisfying these conditions. Then
$$
h_j(T_k x) - \tilde{h}_j(T_k x) = 0 \ \ \mbox{for all $j,k$ and $x$}.
$$
Therefore, the function $h_j - \tilde{h}_j$ is $\Lambda$-periodic and descends to a function on $M$.

If
$$
\psi(x) = e^{i\sum_m \kappa_m h_m(x)}\varphi(x)
$$
is a Bloch eigenfunction, $L\psi=E\psi$, then
$$
\psi(x) = e^{i\sum_m \kappa_m \tilde{h}_m(x)} [e^{i\sum_m \kappa_m (h_m(x) - \tilde{h}_m(x))} \varphi(x)] =
e^{i\sum_m \kappa_m \tilde{h}_m(x)} \tilde{\varphi}(x).
$$
The converse is also true, which implies the coincidence of spectral manifolds constructed from different sets of functions $\{h_j\}$ and $\{\tilde{h}_j\}$.
The theorem has been proven.

{\sc Example.} If $M=S^1$, then $\Gamma/\{\kappa \sim \kappa + 2\pi\}$ is 
the spectral curve of the one-dimensional Schr\"odinger operator with a periodic potential, and the classical Bloch eigenfunctions of this operator \cite{N} correspond to the case $h(x)=x$.

In proof of Theorem 1, we followed our works \cite{T98,T05}. This proof for the periodic Schr\"odinger operator on $\rr^N$ and the heat operator $\partial_y - \partial^2_x + u$, included in the Lax representation for the Kadomtsev--Petviashvili II equation, was obtained by us in 1984 and was even used and cited in \cite{Krichever}, but we published it in \cite{T98} when we found an application of this method for defining the spectral curves of tori immersed in the three-dimensional space.
Note that for the operator $i\partial_y - \partial^2_x +u$, which enters the Lax representation for the Kadomtsev--Petviashvili equation I, the spectral curve does not exist.

\section{Spectra of magnetic Laplacians and the Aharonov--Bohm fluxes}

This work arose as a continuation of \cite{T23} in the study of magnetic Laplacians on non-simply connected manifolds.

The equation \eqref{eqphi1} is rewritten as
$$
(-\widetilde{\Delta} + U)\varphi = E\varphi,
$$
where
$$
\widetilde{\Delta} = \frac{1}{\sqrt{g}}\left(\frac{\partial}{\partial x^j} + i\cA_j\right)\sqrt{g}g^{jk}
\left(\frac{\partial}{\partial x^k} + i\cA_k\right).
$$
The operator $\widetilde{\Delta}$ is the magnetic Laplacian (the Laplace--Boch\-ner operator) corresponding to the magnetic field
$F = \sum_{j<k} F_{jk}dx^j \wedge dx^k$ with potential
$$
\cA = \cA_k dx^k, \ \ d\cA = F.
$$

In our case, for real values of $\kappa_1,\dots,\kappa_N$ we have
$$
\cA = \sum_m \kappa_m \omega_m
$$
and
$$
d\cA = 0,
$$
i.e. this operator corresponds to the magnetic Laplacian with a zero magnetic field.
However, with such a formal procedure for constructing the magnetic Laplacian, the spectrum of $L$ depends on the choice of potential.

The theory of the Aharonov--Bohm effect \cite{AB} also considers the situation with a zero magnetic field on a non-simply connected manifold, which is the plane $\rr^2$ with a certain number $N$ of  punctured points $\{Q_1,\dots,Q_N\} $ (in reality, of small disks). We have
$$
H^1(\rr^2 \setminus\{Q_1,\dots,Q_N\}; \rr) = \rr^N
$$
and we can take the following closed $1$-forms as generators
$\omega_1,\dots,\omega_N$, of this group, such that the integral of $\omega_j$ along a small closed contour $\gamma_k$,
inside which lies only one punctured point --- $Q_k$, is equal to
$$
\oint_{\gamma_k} \omega_j =   \oint_{\gamma_k} A_{j,l}dx^l = \delta_{jk}, \ \ j,k = 1,\dots,N.
$$
When quantizing a zero magnetic field on a plane with punctured points, a choice arises of the ``vector potential'' $\cA = \cA_l dx^l$, which can be cohomologically nontrivial:
$$
[\cA] \neq 0 \in H^1(\rr^2 \setminus\{Q_1,\dots,Q_N\}; \rr).
$$
Expanding it into generators $[\omega_1],\dots,[\omega_N]$,
$$
[\Omega] = \kappa_1 [\omega_1] + \dots + \kappa_N [\omega_N],
$$
we obtain as coefficients the so-called Aharonov--Bohm fluxes
$\kappa_1,\dots,\kappa_N$. Aharonov and Bohm indicated that these fluxes are physically observable, which was later confirmed by physical experiments.

When quantizing a system with the Hamilton function
$$
H(x,p) = g^{jk}(x)p_j p_k +U(x)
$$
in a cohomologically trivial magnetic field $F$, where $[F] = = \in H^2(M;\rr)$, on a non-simply connected manifold $M$ the vector potential of the magnetic field $A$ is expanded into the following sum
$$
\cA = \cA_0 + \sum_m \kappa_m \omega_m, \ \ \ d\cA_0 = F, \ \ d\cA_j = 0, j=1,\dots,N.
$$
It is natural to call the expansion coefficients $\kappa_1,\dots,\kappa_N$ the Aharonov--Bohm fluxes.

Analogously to Theorem 1, the following fact is proved.

\begin{theorem}
When quantizing the motion of a charged particle in an (exact) magnetic field $F$, $dF=0$, on a non-simply connected manifold
spectrum of the corresponding magnetic Laplacian
$$
\widehat{H} = - \frac{1}{\sqrt{g}}\left(\frac{\partial}{\partial x^j} + i\cA_j\right)\sqrt{g}g^{ jk}
\left(\frac{\partial}{\partial x^k} + i\cA_k\right) +U(x)
$$
is related to Aharonov--Bohm flows by the analytical relation
\begin{equation}
\label{ab}
\cP(E,\kappa_1,\dots,\kappa_N) = 0.
\end{equation}
\end{theorem}

Note that for $F \neq 0$ the choice of vector potential $\cA_0$ is ambiguous. The equation
$$
\cP(E,0,\dots,0) = 0
$$
specifies the spectrum of the magnetic Laplacian with $\cA = \cA_0$ and the addition of non-trivial values of the Aharonov--Bohm fluxes shows how the spectrum changes, i.e., in fact, the spectrum at the selected value of the vector potential $\cA_0$ is compared with its deformations when adding to $\cA$ closed $1$-forms. The difference between two vector potentials is a closed $1$-form and the flux values $\kappa_m$ correspond to it.

This relation \eqref{ab} is non-trivial, which is already shown by the example of a zero magnetic field $F=0$.

\begin{theorem}
On an $N$-dimensional torus with $A_0=0, A_{j,l}= \delta_{jl} \kappa_l, j=1,\dots,N$, the equation
\eqref{ab} defines the Bloch variety of the Schr\"o\-din\-ger operator $L = -\Delta+U$ with periodic potential $U(x)$. In this case, the flows $\kappa_1,\dots,\kappa_N$ are quasimoments of Bloch eigenfunctions.
\end{theorem}

According to the Bethe--Sommerfeld conjecture, proven in \cite{Veliev1,Veliev2,Parnovski,Skriganov1,Skriganov2}, for $N \geq 2$ there are a finite number of gaps on the real line such that for any $E$ lying outside these gaps, there is a set of fluxes 
$\kappa_1,\dots,\kappa_N$,
that $E$ is the eigenvalue of the operator
$$
-\sum_{m=1}^N \left(\frac{\partial}{\partial x^m} - i\kappa_m \right)^2 + U(x),
$$
those. Different quantizations of even a zero magnetic field lead to different eigenvalues of the magnetic Laplacian.

\section{Spectral varieties of Dirac operators on Riemann surfaces}

Let us consider the Dirac operator on a two-dimensional spinor Riemannian manifold $M$—a closed Riemann surface of genus $g$.

Let $z$ be a conformal parameter on $M$ and
$$
ds^2 = e^{2\alpha} dz d\bar{z}
$$
be a Riemannian metric. The Dirac operator has the form
$$
D_0 = e^{-3\alpha/2}
\begin{pmatrix} 0 & \partial \\ -\bar{\partial} & 0
\end{pmatrix} e^{\alpha/2},
$$
where $\partial = \frac{\partial}{\partial z}$ and $\bar{\partial} = \frac{\partial}{\partial \bar{z}}$.
It acts on sections of spinor bundles.

We will consider a more general operator by including the potential $V(x)$ --- a matrix function with scalar components that acts on spinors by multiplying:
\begin{equation}
\label{dirac1}
D = D_0 + V = 2 e^{-3\alpha/2}
\begin{pmatrix} 0 & \partial \\ -\bar{\partial} & 0
\end{pmatrix} e^{\alpha/2} + V.
\end{equation}

Since $D_0$ has a discrete spectrum, there is a value of $E_0$ such that the operator
$$
D_0 - E_0,
$$
acting on $L_2$-sections of the spinor bundle is invertible. We can apply the same reasoning as in the proof of Theorem 1 to derive the following result.
The only serious difference is the analogue of Kondrashov's theorem on compactness of embedding for sections of spinor bundles, but in the case of compact manifolds it immediately follows from a similar theorem for function spaces.

\begin{theorem}
In $\cc^{2g+1}$ there is an analytic subset, the spectral variety, determined as the zero set
$$
S = \{\cP(E,\kappa_1,\dots,\kappa_N)=0\}
$$
of a nontrivial entire function such that $(E,\kappa_1,\dots,\kappa_{2g})$ lies to
$S$ if and only if there is a Bloch eigenspinor of the form 
\begin{equation}
\label{bs}
\psi(x) = e^{i(\kappa_1 h_1(x) + \dots + \kappa_{2g} h_{2g}(x))} \varphi(x),
\end{equation}
on $M_0$, 
where $\varphi$ is a section of the spinor bundle over $M$,
satisfying the equation
$$
D\varphi = E\varphi.
$$

The spectral variety does not depend on the choice of functions $h_1$, $\dots$, $h_{2g}$ such that
$p^\ast \omega_k = dh_k, k=1,\dots,2g$, and therefore is determined only by the operator $D$.
\end{theorem}

\section{Spectral varieties as invariants of surfaces in Euclidean three- and four-spaces}

\subsection{The Weierstrass (spinor) representation of surfaces in three- and four-spaces}

Any immersion of a closed oriented surface $M$ immersed into $\rr^3$ or $\rr^4$ is described in terms of a solution to the Dirac equation
$$
\cD\psi = 0,
$$
where
$$
\cD = \begin{pmatrix} 0 & \partial \\ -\bar{\partial} & 0 \end{pmatrix} + \begin{pmatrix} U & 0 \\ 0 & \bar{U} \end{pmatrix}.
$$
For brevity, we will limit ourselves to surfaces in $\rr^3$. In this case, the potential $U$ is real. Following \cite{T15}, we identify $\rr^3$ with the space of Hermitian $2 \times 2$-matrices $A$: 
$\bar{A}^\top = -A$, and represent the immersion in the form
\begin{equation}
\label{weier}
\begin{split}
X(P) = \begin{pmatrix} ix^3 & -x^1 - ix^2 \\ x^1-ix^2 & -ix^3 \end{pmatrix} (P) =
\\
=
i \int_{P_0}^P \left[
\begin{pmatrix} \psi_1 \bar{\psi}_2 & -\bar{\psi}_2^2 \\ \psi_1^2 & -\psi_1 \bar{\psi}_2 \end{pmatrix} dz +
\begin{pmatrix} \bar{\psi}_1 \psi_2 & \bar{\psi}_1^2 \\ -\psi_2^2 & -\bar{\psi}_1 \psi_2
\end{pmatrix} d\bar{z} \right] + X(P_0).
\end{split}
\end{equation}
Here $x^1,x^2,x^3$ are the Euclidean coordinates in $\rr^3$, $\psi = \begin{pmatrix} \psi_1 \\ \psi_2 \end{pmatrix}$ is a solution of the Dirac equation, $X: \cU \to \rr^3$ is an immersion of a simply connected domain $\cU \subset \cc$ in $\rr^3$, $z$ is a complex parameter (coordinate) in $\cU$, $P_0$ is a fixed point in $\cU$ and the integral is taken along any path from $P_0$ to $P$ lying in 
$\cU$ (its value does not depend on the path).

The vector function $\psi$ defines the surface with up toa shift $X(P_0)$, i.e. 
essentially defines its Gaussian map.

These formulas for the local representation of surfaces allowing soliton deformations were derived in \cite{Kon}, where it was shown that the parameter $z$ is conformal, the induced metric takes the form
$$
ds^2 = (|\psi_1|^2 + |\psi_2|^2)^2 \,dz\,d\bar{z} = e^{2\alpha}\,dz\,d\bar{z }
$$
and the potential $U$ of the Dirac operator takes the form
$$
U = \frac{e^\alpha}{2} H,
$$
where $H$ is the mean curvature of the surface (in the case of surfaces in $\rr^3$, the potential $U$ is real).

In \cite{T97} it was proven that any closed oriented surface in $\rr^3$ has such a representation, and that $\psi_1$ and $\psi_2$ are sections of a spinor bundle over the immersed surface.

\subsection{Spectral curves of immersed tori.} 
Let $M$ be a torus immersed in $\rr^3$ by \eqref{weier}. Let
$$
\cc/\{\zz e_1+ \zz e_2\} \approx M
$$
be a torus conformally equivalent to $M$. We denote by $e_1$ and $e_2$ the generators of the period lattice $\Lambda$. We will identify them with the generators of $H_1(M;\zz_2)$.

The immersion is given by the spinors $\psi_1 \sqrt{dz}$ and $\bar{\psi}_2\sqrt{dz}$. The spinor bundle over any closed oriented Riemann surface is given by the quadratic form
$$
q: H_1(M;\zz_2) \to \zz_2,
$$
where, briefly speaking, $p(w) = (-1)^{q(w)}$ is the monodromy of a spinor under a shift by $w$.
Recall that a form on $H_1(M;\zz_2)$ is called quadratic if
$$
q(w_1 + w_2) = q(w_1) + q(w_2) + w_1 \cdot w_2,
$$
where $w_1 \cdot w_2$ is the intersection (modulo two) of the cycles $w_1$ and $w_2$.

For example, it is easy to show that $\sqrt{dz}$ is a section of the bundle for which
$$
q_0(0) = 0, \ \ q_0(e_1) = q_0(e_2) = q_0(e_1+e_2) = 1.
$$

Since the quantities $\psi_1^2 dz, \bar{\psi}_2^2 dz$ and $\psi_1 \bar{\psi}_2 dz$ are doubly periodic, then
$\psi$ is the Bloch eigenfunction of $\cD$, which, when shifted by periods, is multiplied by
$\mu(w)=\pm 1$. Note that the multiplier mapping is given by the homomorphism
$$
\nu: H_1(M;\zz_2) \to \zz_2, \ \ \ \mu(w) = (-1)^{\nu (w)}.
$$
The spinors $\psi_1 \sqrt{dz}$ and $\bar{\psi}_2 \sqrt{dz}$ are sections of the spinor bundle,
which corresponds to the quadratic form
$$
q = q_0 + \nu.
$$
There are four different homomorphisms $H_1(M;\zz_2) \to \zz_2$, which, using the previous formula, define four different quadratic forms on $H_1(M;\zz_2)$ (and thereby spinor structures on $M$) .

From the above it follows that the spectral variety is determined independently of the spinor structure. 

In \cite{T98} we introduced the spectral curve of an immersed torus as the Bloch variety of $\cD$ with doubly periodic potential $U(x)$ and at zero energy level $E=0$. Bloch eigenfunc\-tions are solutions to the equations
$$
\cD \psi = 0, \ \
\psi (x + e_j) = e^{2\pi i \langle k, e_j \rangle} \psi(x) = \mu(e_j) \psi(x), \ \ j=1,2.
$$
After replacement
$$
\psi \to \varphi = e^{-\alpha/2}\psi
$$
this Dirac equation reduces to
$$
D\varphi = 0,
$$
where the operator $D$ has the form \eqref{dirac1} and
$$
V = \begin{pmatrix} 2U e^{-\alpha/2} & 0 \\ 0 & 2 U e^{-\alpha/2} \end{pmatrix} =
\begin{pmatrix} H & 0 \\ 0 & H \end{pmatrix}.
$$
It follows from Theorem 4 (after the obvious linear transformation $\kappa \to k$) that
such solutions exist, exactly, for $\cP(k_1,k_2) = 0$, where $\cP$ is some entire function.
Since the multipliers $\mu$ are invariant under the shift $k$ by vectors from the dual lattice:
$$
\Lambda^\ast = \{ v \ : \langle v, e_k\rangle \in \zz, k=1,2 \},
$$
then after factorization
$$
S = \{\cP=0\} \to \{ \cP=0\} / \Lambda = \Gamma
$$
we obtain a Riemann surface $\Gamma$, which parametrizes all Bloch eigenfunc\-tions of the operator $\cD$.

This spectral curve is ``almost'' preserved under conformal transformations of the ambient space. Namely, the  multipliers are preserved.
The simplest proof is the following: it is quite obvious that the multipliers are preserved when $\rr^3$ rotates (they reduce to the actions of $SO(3)$ on $\psi$). The inversion in terms of \eqref{weier} takes the form
$$
X \to X^{-1}
$$
and there is a Moutard transformation connecting the Weierstrass representations of tori and their Bloch eigenfunctions while preserving the multipliers \cite{T15}. At the same time, \cite{MT} contains examples where inversion leads to the appearance of double points on the spectral curve.

For very large quasimomenta, the spectral curve ``approaches'' the spectral curve in the case of zero potential:
$$
k_2 \approx \pm i k_1.
$$
Near the asymptotic end, where $k_2 \approx ik_1$ we can introduce such a local parameter
$\lambda^{-1}$ that
\begin{equation}
\label{will1}
\mu(v) = \lambda v + \frac{C_0 \bar{v}}{\lambda} + O(\lambda^{-2}),
\end{equation}
where
\begin{equation}
\label{will2}
C_0 = -\frac{1}{\mathrm{Area}\, (\cc/\Lambda)} \int_{\cc/\Lambda} U^2 =
-\frac{1}{4\, \mathrm{Area}\, (\cc/\Lambda)} \int_M H^2 d\mu
\end{equation}
is $\int_M H^2 d\mu$ is the Willmore functional of $M$.

This observation formed the basis of our proposed ``spectral'' approach to Willmore's conjecture, which remained unrealized (see discussion in \cite{T97,T05,T23}).

Note that the spectral curves of immersed tori can be reducible \cite{BT2}.

All these statements are also true for surfaces in $\rr^4$. Only in this case can the potential $U(x)$ be complex-valued and the Weierstrass representation is given by a pair of spinors $\psi$ and $\phi$ such that $\cD\psi = 0$ and $\cD^\vee \phi= 0$, where the operator $\cD^\vee$ is obtained from $\cD$ by permuting $U$ and $\bar{U}$ on the diagonal in the matrix coefficient (see, for example, \cite{T05,T06}). The Moutard transformation corresponding to inversions in $\rr^4$ was derived in \cite{MT}.

\subsection{Spectral varieties of higher genera surfaces.}
The arguments from \S 5.2 extend to surfaces of large genus immersed in $\rr^3$ and $\rr^4$.
We have

\begin{theorem}
If $D$ is a Dirac operator of the form \eqref{dirac1}, included in the Weierstrass representation
closed oriented surface of genus $g$ immersed into $\rr^3$ or $\rr^4$, then
analytic set
$$
\{\cP(0,\kappa_1,\dots,\kappa_{2g}) = 0\},
$$
defined by Theorem 4 and factorized by the action $\kappa \to \kappa + 2\pi n, n \in \zz$,
is the spectral variety of $\Gamma$. 
The multiplier function $\mu =
(e^{i \kappa_1},\dots,e^{i\kappa_{2g}})$ invariant under conformal transformations of the ambient space.
\end{theorem}

Probably these complex manifolds should contain information about the Willmore functional of the surface similar to that given in \eqref{will1} and \eqref{will2}. However, already for $g=2$ the spectral variety has a complex dimension equal to three. If a surface of large genus has symmetries, then this should lead to special classes of spectral varieties.

\vskip7mm
The work was carried out with the support of the Mathematical Center in Akademgorodok,
agreement with the Ministry of Science and Higher Education of the Russian
Federation No. 075-15-2022-281.

\end{document}